\documentclass[a4paper]{amsart}
\usepackage[all]{xy}

\newtheorem{num}{}[section]
\theoremstyle{definition}
\newcommand{\rmap}{\longrightarrow}

\begin{document}

\title [On the perturbation lemma]{On the perturbation lemma, and deformations}

\author{Marius Crainic}
\address{Depart. of Math., Utrecht University, 3508 TA Utrecht, 
The Netherlands}
\email{crainic@math.uu.nl}

\date{}
\begin{abstract} 
We have one more look at the perturbation lemma and we point out some 
non-standard consequences, including the relevance to deformations.
\end{abstract}
\maketitle

\section{Introduction}

The classical (homological) perturbation lemma is a technical tool from homological algebra which is very useful both ``for proving'', as well as ``for providing explicit formulas''. It can be viewed as an algebraic (or homological) version of Newton's iteration method. For several applications of the classical perturbation
lemma, as well as for more references, please see \cite{G, HuK, Kas, St}. For a short presentation, see also 
Remark \ref{in-rem} (iii) below. Here we improve the classical statement and we point out some consequences, which, previously, did not appear to be related to the perturbation lemma. In relation with deformations, the results derived here are just an indication of the relevance of the lemma (and its variations) to deformation theory (see also Remark \ref{def-Lie-rk} and Remark \ref{def-Ban-rk}).

{\it Acknowledgments:} I would like to thank R.L. Fernandes for his comments on a first version of this paper. 

\section{The main perturbation lemma}
\label{TMPL}

\begin{num}{\bf Initial data:}\rm \  A {\it HE data} (homotopy equivalence data):

\begin{equation}
\label{unu}  \xymatrix{
(L, b ) \ar@<-1ex>[r]_-{i} &   (M, b),\ h \ar[l]_-{p}}
\end{equation}
consists on the following:

\begin{enumerate}
\item[(i)] two complexes $(L, b )$, $(M, b)$ and q.i.'s $i$, $p$ between them;
\item[(ii)] a homotopy $h$ between $ip$ and $1$ (so $ip= 1+ bh+ hb$);
\end{enumerate}
By complexes we mean here either chain or cochain complexes. Recall also that ``q.i.'s'' stands for ``quasi-isomorphisms'', i.e. (co)chain maps which induce isomorphism in (co)homology. Unless the contrary is specified, all maps are assumed to be linear.
\end{num}

\begin{num}{\bf The perturbed data:}\rm \  A {\it perturbation} $\delta$ of (\ref{unu}) is a map on $M$ of the same degree as $b$,  such that $(b +\delta)^{2}= 0$.
We call it {\it small} if $(1- \delta h)$ is invertible. In this case we put:
\[ A= (1- \delta h)^{-1} \delta ,\]
and we consider
\begin{equation}
\label{trei} 
\xymatrix{
(L, b_{1}) \ar@<-1ex>[r]_-{i_{1}} &   (M, b+ \delta ),\ h_{1} \ar[l]_-{^{p_{1}}} }
\end{equation}
with:
\[ i_{1}= i+ hAi, \ \ \  p_{1}= p+ pAh,\ \ \ h_{1}= h+ hAh, \ \ \  b_{1} = b+ pAi. \]
\end{num}

\begin{num}\label{in-rem}{\bf Remarks:}\rm \ 
\begin{enumerate}
\item[(i)] A DR (deformation retract) is a HE data (\ref{unu})
with the property that $pi= 1$. Note that this relation does not survive after perturbation.
To fix this, one can restrict himself to {\it special DR}'s, i.e. to those DR which also satisfy
\[ hi= 0, ph= 0, h^{2}= 0 .\]
Any DR can be made into a special DR by the following transformations: replace $h$ by 
$h(bh+hb)$ to realize $hi= 0$, replace $h$ by $(bh+hb)h$ to realize $ph= 0$, and then replace $h$ by
$hbh$ to realize $h^2= 0$. Note however that, when interested on explicit formulas, such transformations are often un-wanted: they do affect the perturbed data,
and the formulas become much more complicated.
\item[(ii)] One can check that HE's enjoy the same functorial properties as special DR's \cite{St}. Note also that $-\delta$ can be regarded
as a small perturbation of the perturbed data, and, perturbing once more, one recovers the initial data.
\item[(iii)] The classical perturbation lemma says that the perturbed data
is a special DR provided the initial data is a special DR, and there exists
a filtration of (\ref{unu}) (i.e. filtrations of $L$ and $M$ which are preserved by $i, p$ and $h$),
which is increasing and bounded below, and which is lowered by $\delta$. These conditions imply that 
$\delta h$ is locally nilpotent, i.e., for each $x\in M$, there exists  $n\geq 0$ such that
 $(\delta h)^n(x) = 0$. In turn, this implies that $\sum_{n= 0}^{\infty} (\delta h)^n$ is defined
proving that $\delta$ is small. Other variations of
the perturbation lemma are based on the same type of conditions.
\item[(iv)] Note that, with the mind at the formal equation
\[ (1- \delta h)^{-1}=  \sum_{n= 0}^{\infty} (\delta h)^{n} ,\]
there are other various particular cases when $\delta$ is small. In particular, in the presence of Banach norms,
small enough $\delta$'s are small in the previous sense. 
\end{enumerate}
\end{num}

\begin{num}
\label{LLL}
\underline{{\bf Main Perturbation Lemma:}} If $\delta$ is a small perturbation of the HE data (\ref{unu}), then
the perturbed data (\ref{trei}) is a HE.
\end{num}

Let (\ref{unu}) be a HE and $\delta$ a perturbation as in the statement. We need the following:

{\bf Lemma:} {\it We have the following relations}:
\begin{eqnarray}
\lefteqn{\hspace{-1.3in} \delta hA= Ah\delta= A- \delta, \label{i}}\\
\lefteqn{\hspace{-1.3in} (1-\delta h)^{-1}= 1+ Ah,\ \ \  (1- h \delta)^{-1}= 1+ hA, \label{ii}}\\
\lefteqn{\hspace{-1.3in} AipA+ Ab+ bA= 0. \label{iii}}
\end{eqnarray}

\begin{proof} From the definition of A, $(1- \delta h)A= \delta$ which proves $\delta hA= A-\delta$. Multiplying the identity $\delta h\delta =\delta - (1- \delta h)\delta$ by $(1- \delta h)^{-1}$ from the left we also get $Ah\delta =A-\delta$. These prove (\ref{i}). The relations we have to check in order to prove  (\ref{ii}) follows immediately from (\ref{i}). For instance:
\[ (1-\delta h)(1+ Ah)= 1+ Ah- \delta h- \delta hAh= 1+ h(A- \delta- \delta hA)= 1.\]
To prove (\ref{iii}) we use (\ref{i}), (\ref{ii}) and the relations $ip= 1+ bh+ hb, \delta^{2}+ \delta b+ b\delta =0:$
\begin{eqnarray*}
AipA+ bA+ Ab & = & A(1+ bh+ hb)A+ bA+ Ab\\
             & = & A^{2}+ Ab(hA+ 1)+ (Ah+1)bA  \\
             & \stackrel{  (\ref{ii})}{=} & A^{2}+ Ab(1- h\delta)^{-1}+ (1- \delta h)^{-1}bA \\
             & = & (1-\delta h)^{-1}[ (1-\delta h)A^{2}(1- h\delta)+ (1- \delta h)Ab+ bA(1- h\delta) ](1-h\delta)^{-1}\\
             & = & (1-\delta h)^{-1}[ (A- \delta hA)(A- Ah\delta)+ (A-\delta hA)b+ b(A- Ah\delta) ](1-h\delta)^{-1}\\
             & = & (1-\delta h)^{-1}[ \delta \delta + \delta b + b \delta ](1-h\delta)^{-1} = 0 . \ \   
\end{eqnarray*}
\end{proof}

\begin{proof} ({\it of the Main Perturbation Lemma})  We have to prove various relations:\\ \newline
\hspace*{.3in} $\underline{1). b_{1}^2= 0}$\ \ (i.e. $b_{1}$ is, indeed, a boundary):
\begin{eqnarray*}
b_{1}^2 & = & (b+ pAi)(b+ pAi) = b^2+ bpAi+ pAib+ p(AipA)i\\
             &\stackrel{  ( \ref{iii})}{=}&  bpAi+ pAib- p(Ab+ bA)i= 0.
\end{eqnarray*}
\hspace{.3in} $\underline{2). i_{1}b_{1}= (b+ \delta) i_{1}}$\ \ (i.e. $i_{1}$ is, indeed, a chain map):
\begin{eqnarray*}
i_{1}b_{1}- (b+ \delta) i_{1} & = & (i+ hAi)(b+ pAi)- (b+\delta)(i+ hAi)\\
  & = &  ib+ ipAi+ hAib+ h(AipA)i- bi- bhAi- \delta i- (\delta hA)i\\
  & \stackrel{  ( \ref{i}),   ( \ref{iii})}{=} & ipAi+ hAib- h(Ab+ bA)i- bhAi- \delta i- (A-\delta)i\\
  & = & ipAi- hbAi- bhAi- Ai = (ip- hb- bh- 1)Ai = 0 .
\end{eqnarray*}
\hspace{.3in} $\underline{3). b_{1}p_{1} = p_{1}(b+ \delta)}$\ \ (i.e. $p_{1}$ is, indeed, a chain map):
\begin{eqnarray*}
b_{1}p_{1} - p_{1}(b+ \delta) & = & (b+ pAi)(p+ pAh)- (p+ pAh)(b+ \delta)\\
\hspace{.3in}   & = & bp+ bpAh+ pAip+ p(AipA)h- pb- p\delta - pAhb- p(Ah\delta)\\
    & \stackrel{  ( \ref{i}),   ( \ref{iii})}{=}  & bpAh+ pAip- p(Ab+ bA)h- p\delta - pAhb- p(A-\delta)\\
    & = & pAip- pAbh- pAhb -pA = pA(ip- bh- hb- 1) = 0
\end{eqnarray*}
\hspace{.3in} $\underline{4). i_{1}p_{1} = 1+ h_{1}(b+ \delta)+ (b+ \delta)h_{1}}$\ \ (i.e. $h_{1}$ is a homotopy between $i_{1}p_{1}$ and $1$):
\begin{eqnarray*}
\lefteqn{\hspace*{.5in} 1+ h_{1}(b+ \delta)+ (b+ \delta)h_{1}- i_{1}p_{1} = }\\
& = & 1+ (h+hAh)(b+\delta)+ (b+\delta)(h+ hAh)- (i+ hAi)(p+ pAh)\\
& = & \underline{1+ hb}+ h\delta+ hAhb+ h(Ah\delta)+ \underline{bh}+ \delta h+ bhAb+ (\delta hA)h- \underline{ip}- ipAh- hAip- h(AipA)h\\
& \stackrel{  ( \ref{i}),   ( \ref{iii})}{=} & h\delta + hAhb+ h(A-\delta)+ \delta h+ bhAh+ (A-\delta)h- ipAh- hAip+ h(Ab+ bA)h\\
& = & \underline{hAhb+ hA}+ bhAh+ Ah- ipAh- \underline{hAip}+ \underline{hAbh}+ hbAh\\
& = & hA(hb+ 1- ip+ bh)+ (bh+ 1- ip+ hb)Ah = 0.
\end{eqnarray*}
\hspace*{.3in} $\underline{6). p_{1} \ {\rm and}\  i_{1} \ {\rm are\ quasi-isomorphisms}}$: From step 4 it follows that $i_1p_1$ induces the identity in homology. So it suffices to show that $i_1$ is injective in homology. Assume that $x\in L$ has $b_1(x)= 0$ and $i_1(x)= (b+\delta)y$ for some $y\in M$. Hence
\begin{eqnarray}
 & & b(x)+ pAi(x)= 0 \label{e1}\\
 & & i(x)+ hAi(x)= b(y)+ \delta(y) \label{e2}
\end{eqnarray}
Applying $\delta$ to the last equation, and replacing $\delta hA$ by $A-\delta$ (cf. Lemma) and $\delta b+ \delta^2$ by $-b\delta$, we obtain
\[ Ai(x)= -b\delta(y) .\]
With this formula for $Ai(x)$ plugged into (\ref{e2}), we get 
\begin{equation}
\label{**}
i(x)= b(y)+ \delta(y)+ hb\delta(y),
\end{equation} and, using
$hb= ip- 1- bh$, we get 
\[ i(x- p\delta(y))= b(y- h\delta(y)).\]
Next, plug the formula (\ref{*}) for $Ai(x)$ into (\ref{e1}) to get 
\[ b(x- p\delta(y))= 0.\] 
The last two formulas tell us that $x- p\delta(y)$ is a cocycle which, in cohomology, is mapped by $i$ into zero. 
Hence, since $i$ is a q.i., we find $z\in M$ such that
\begin{equation}
\label{*}
x= p\delta(y)+ b(z) .
\end{equation} 
Applying $i$ to this formula and using (\ref{**}),
\[ b(y)+ \delta(y)+ hb\delta(y)= ip\delta(y)+ ib(z).\]
Using now $ip= 1+ hb+ bh$, we deduce that
\[ b(i(z)- y+ h\delta (y))= 0.\]
Since $i$ is surjective in cohomology, we can write
\[  i(z)- (1- h\delta)(y)= i(\alpha)+ b(\beta)\]
for some $\alpha\in L$ with $b(\alpha)= 0$, and some $\beta\in M$. Applying $pA$ to this and using $A(1- h\delta)= 1$ and 
$Ab= -bA- AipA$ (cf. Lemma), we deduce
\[ pAi(z)= p\delta(y)+ pAi(\alpha)- pbA(\beta)- pAipA(\beta) .\]
From this we extract $p\delta(y)$ and we plug the result in (\ref{*}). Rearranging the terms we get
\[ x= b(z+ pA(\beta))+ pAi(z+ pA(\beta)-\alpha)) .\]
Since $b(\alpha)= 0$, we conclude that $x= b_1(z+ pA(\beta)-\alpha)$ is a exact.
\end{proof}

\section{Immediate consequences and variations}

In this section we point out some consequences and variations. We start with some very simple consequences
(\ref{Perturbing DR's}- \ref{main-cor}) with the aim of showing the 
usefulness of the main perturbation lemma as ``a very general recipe'' for obtaining explicit formulas
(hence the use of the lemma will hopefully lower our use of aspirin).
Then we have a look at topological versions.

\begin{num}\label{Perturbing DR's} {\bf Perturbing DR's}:
If $\delta$ is a small perturbation of the DR data (\ref{unu}), then the perturbed data 
is a DR if and only if $p\, (Ah^{2}A+ Ah +hA\, ) i= 0$.
\end{num}

For special DR's, we obtain:

\begin{num} {\bf Perturbing special DR's}:
If (\ref{unu}) is a special DR, and $\delta$ is a perturbation of $b$ so that $(1-\delta h)$ is invertible,
then (\ref{trei}) is a special DR.
\end{num}

The next particular case (when $h= 0$), although trivial (or especially because of that), is a
very good illustration of one type of examples that the classical perturbation lemma does not handle (due 
to the requirement that $pi= 1$).

\begin{num}\label{retr} {\bf Corollary}: Let $\delta$ be a perturbation of the cochain complex $(M, b)$ and let 
\begin{equation}
\xymatrix{
(L, b ) \ar@<-1ex>[r]_-{i} &   (M, b) \ar[l]_-{p}}
\end{equation}
with $i$ and $p$ q.i.'s satisfying $ip= 1$. Then 
\begin{equation}
\xymatrix{
(L, b+ p\delta i ) \ar@<-1ex>[r]_-{i} &   (M, b+ \delta) \ar[l]_-{p}}
\end{equation}
has the same properties.
\end{num}

Another interesting case is when $L= 0$:

\begin{num}\label{contr} {\bf Perturbing contractions}:
If $(M, b)$ is a contractible complex with contraction $h$ (i.e. $hb+bh+1= 0$),
and if $\delta$ is a perturbation of $b$ (i.e. $(b+\delta)^2= 0$) such that $(1- \delta h)$ is invertible, then $(M, b+\delta)$
is contractible, with contraction 
\[ H= h(1-\delta h)^{-1}.\]
\end{num}

The following is probably part of the folklore, but it is often very useful for finding explicit formulas,
e.g. for the Chern character in the various complexes computing the cyclic cohomology (where the obvious DR's
to consider are not special). 

\begin{num}{\bf Contracting rows}: 
Let $C_{p, q} \ (p, q \geq 0)$ be a double complex with horizontal differentials $\partial$  and vertical differentials $\delta$ ($\partial^2= \delta^{2}= \partial \delta+ \delta \partial = 0$). Assume that each row $C_{*, q}$ has zero cohomology in positive degrees, and that
we choose contractions of the augmented rows:
\[ \xymatrix{
H_q= H_{0}(C_{*, q}) \ar@<-1ex>[r]_-{i} & C_{0,q} \ar@<-1ex>[r]_-{h} \ar[l]_-{p} & C_{1,q} \ar@<-1ex>[r]_-{h} \ar[l]_-{\partial} & C_{2,q} \ar[l]_-{\partial} \ar@<-1ex>[r]_-{h}& .\ .\ .\ \ar[l]_{\partial} .}\]
(i.e. $pi= 1$, $ip= \partial h+ h\partial + 1$). Then:
\[ \xymatrix{ (H_{*}, \delta) \ar@<-1ex>[r]_-{i_{1}} & ( Tot(C), \partial + \delta), \ h_{1} \ar[l]_-{p} } \]
is a DR, where:
\[ i_{1}= \sum_{n \geq 0}   ( h\delta )^{n}i, \ \ h_{1}= \sum_{n \geq 0} h  ( \delta h)^{n} .\]
\end{num}

\begin{proof} Apply \ref{LLL} to $(H_{*}, 0), (Tot(C_{*, *}, \partial)$ with $\delta$ viewed as perturbation, with the obvious projection $p$ and inclusion $i$, and with the homotopy $h$. The formula $p(Ah^2A+ Ah+ hA)i= 0$ is automatically satisfied.  
\end{proof}

The following is a basic lemma in constructing the $b-B$ complex computing the cyclic homology (Lemma 2.1.6 in \cite{Lo}), and is another example
where the role of the perturbation lemma is to produce explicit formulas. 

\begin{num} {\bf Killing Contractible Complexes}:
Let
\[ .\ .\ . \rmap A_{n}\oplus A_{n}' \stackrel{d}{\rmap} A_{n-1}\oplus A_{n-1}' \rmap .\ .\ . \ ; d= \left[ \begin{array}{cc}
                                               \alpha & \beta \\
                                               \gamma & \delta
                                                   \end{array}   \right]  \]
be a chain complex such that $(A'_{*},\delta)$ is contractible with contracting homotopy $H: A_{n}' \rmap A_{n+1}'$ ($H\delta +\delta H= 1$). Then the following is a HE:
\[  \xymatrix{ (A_{*}, \alpha- \beta H\gamma) \ar@<-1ex>[r]_-{i_{1}} & (A_{*}\oplus A'_{*}, d), \ h_{1} \ar[l]_-{^{p_{1}}} }, \]
\[ i_{1}  ( a\, ) = (a, - H\gamma  ( a\, )),\ \  p_{1}  ( a, b\, )= a- \beta H  ( b\, ) ,\ \ h_{1}= \left[ \begin{array}{cr}
                           0 & 0\\
                           0 & -H
                     \end{array} \right]. \]
\end{num}

\begin{proof} Use \ref{LLL} for $(A_{*}, 0), (A_{*}\oplus A_{*}', b)$, $b= \left[ \begin{array}{cc}
                                               0 & 0 \\
                                               0 & \delta
                                                   \end{array}   \right] $ with the obvious projection and inclusion, the homotopy $h= \left[ \begin{array}{cr}
                           0 & 0\\
                           0 & -H
                     \end{array} \right]$ and the perturbation $\left[ \begin{array}{cc}
                                               \alpha & \beta \\
                                               \gamma & 0
                                                   \end{array}   \right].\ \   $ \\
\end{proof}

Concentrating on certain degrees only, or looking at topological complexes,  
there are several variations of the main perturbation lemma and of its consequences. 
For instance, let us state the following variation of \ref{contr}, 
which follows directly from the formulas in the main proof.

\begin{num}\label{main-cor} {\bf Corollary}: Let $(C, b)$ be a cochain complex, and let $h$ be a contraction in degree $k$, i.e.
\begin{equation}  
\xymatrix{
C^{k-1} \ar@<-1ex>[r]_-{b} &   C^k \ar[l]_-{h}\ar@<-1ex>[r]_-{b} & C^{k+1}\ar[l]_-{h} ,\ bh+ hb+ 1= 0\ \text{on}\ C^k.}
\end{equation}
Then, for any perturbation $\delta$ of $b$ with the property that $1- \delta h$ is invertible on $C^k$ and $C^{k+1}$,  
$h= h(1-\delta h)^{-1}$
define a contraction in degree $k$ of $C$ with the differential $b_1= b+\delta$.
\end{num}

It is well known that the existence of a contraction $h$ in degree $k$ is equivalent to the vanishing of $H^k(C, b)$. 
Indeed, choosing a right inverse $\lambda$ for $b: C^k\rmap Im(b)$, a right inverse $\lambda'$ for $b: C^{k-1}\rmap Im(b)$,
and a left inverse $\pi$ for the inclusion $Im(b)\rmap C^{k+1}$,
\begin{equation}
\label{formula-alg-h} 
h= \left \{ \begin{array}{ll} 
                                              \lambda\pi & \mbox{on $C^{k+1}$} \\
                                              \lambda'(1-\lambda b) & \mbox{on $C^{k}$} 
                                                        \end{array}
                                            \right. 
\end{equation}
defines the desired contraction.

In the topological category however, such an argument is no longer possible, since continuous surjections need not have
right inverses which are both linear and continuous. Hence, the obvious topological version of the previous corollary has to 
assume the existence of $h$. The following is yet another variation, and this one can handle cases where only the cohomology vanishing 
is assumed. Given two l.c.s.'s (locally convex topological vector spaces) $V$ and $W$, a map $h: V\rmap W$ is called quasi-bounded if, 
for any semi-norm $q$ of $W$, there exists $C>0$ and a seminorm $p$ of $V$ such that
\[ q(h(v))\leq C p(v) \]
for all $v\in V$. We have the following

\begin{num}\label{contr-top}{\bf Lemma}:\ Let $(C, b)$ be a cochain complex of l.c.s.'s, and let $h$ be a contraction in degree $k$, where, this time,
$h$ may be non-linear, but we assume that $h$ is quasi-bounded on $C^k$ and $h(-v)= -h(v)$.

If $\delta$ is a perturbation of $(C, b)$ which is ``small enough'' in the sense that  
\[ \sum [\delta, h]^n(v)   \]
is absolutely convergent for all $v\in Ker(b+\delta)\cap C^k$ ,
then $H^k(C, b+\delta)= 0$. 
\end{num}

Here, $[\delta, h]= \delta h+ h\delta$, and a sum $\sum v_n$ in a l.c.s. $V$ is absolutely convergent if
$\sum p(v_n)$ is convergent for all seminorms $p$ of $V$.  

\begin{proof} Let $v\in C^k$ such that $(b+\delta)v= 0$. We show that
\[ v= (b+\delta)(w), \ \mbox{where}\ w= -\sum_{n= 0}^{\infty} h([\delta, h]^{n}(v)).\]
(Recall that $h$ does not commute with sums!). 
First of all, remark that $[\delta, h]$ preserves $Ker(b+\delta)$. Indeed, 
\[ (b+\delta)[\delta, h]= (b+\delta)\delta h+ (b+\delta) h\delta= -\delta b h+ (b+\delta) h\delta,\]
hence, on $Ker(b+\delta)$, this is
\[ -\delta b h- (b+\delta) hb= -\delta (bh+hb)- bhb= \delta+ b= 0.\]
Hence each $v_n:= [\delta, h]^n(v)$ is in $Ker(b+\delta)$. Next,
since $bh= -1+ h\delta$ on $Ker(b+\delta)$, 
\[ (b+\delta) h(v_n)= -v_n+ h\delta(v_n)+ \delta h(v_n)= v_{n+1}- v_{n} .\]
Hence 
\[ (b+\delta) \sum_{i= 0}^{n} h(v_i)= v_{n+1}- v .\]
Since $\sum v_n$ is absolutely convergent it follows that $v_n$ converges to zero.
Since $h$ is quasi-bounded, it also follows that $\sum h(v_n)$ is absolutely convergent. Hence, 
taking $n\rightarrow \infty$ in the last equation, we get $(b+\delta)(w)= v$. 
\end{proof}

%

In particular, we deduce the following (compare with Section 6 of \cite{Jo}).

\begin{num}\label{cor-pert-Banach}{\bf Corollary}: If $(C, b)$ is a cochain complex of Banach spaces so that
$H^k(C, b)= 0$ and $Im(b)\subset C^{k+1}$ is closed, then there exists $\epsilon > 0$ such that, for all perturbations
$\delta$ of $(C, b)$ with $||\delta||< \epsilon$ on $C^{k-1}$ and on $C^k$, $H^k(C, b+\delta)= 0$.
\end{num}

\begin{proof} We proceed as in the algebraic case to construct $h$ by formula (\ref{formula-alg-h}). This time,
$\lambda$ and $\lambda'$ will be chosen to be non-linear but quasi-bounded (this is possible
by the open mapping theorem and the fact that $b$ has closed range in the necessary degrees), and $\pi$ will
still be linear. Also, eventually replacing $\lambda(v)$ by $\frac{1}{2}(\lambda(v)- \lambda(-v))$, we may assume that 
$\lambda(-v)= -\lambda(v)$, and similarly for $\lambda'$. We only have to show that the sum appearing in the
previous statement is absolutely convergent. On the other hand, for $v\in Ker(b+\delta)$, since $\delta(v)\in Im(b)$,
\[ [\delta, h](v)= \delta h(v)+ \lambda \delta(v) ,\]
and, on norm, this is less than
\[ \epsilon C_h ||v||+ C_{\lambda}||\delta(v)||\leq \epsilon (C_h+ C_{\lambda})||v|| ,\]
for some constants $C_h$ and $C_{\lambda}$. Hence, it suffices to choose $\epsilon< 1/(C_h+ C_{\lambda})$.
\end{proof}

\section{Deforming Lie algebras}

From this section on, we look at deformations of algebraic/topological structures. Such deformations are indexed by the real parameter 
$t$ varying in an open interval containing the origin, will depend smoothly on $t$, and, at $t= 0$, one recovers the original structure that is being deformed. Two such deformations are equivalent if, for $t$ small enough, they are isomorphic by isomorphisms $h_t$ depending smoothly on $t$, with $h_0= Id$.  

In this section we illustrate the use of the perturbation lemma to deformations by looking at one of the simplest examples: Lie algebras. This should be viewed as a ``baby example'', but which is suggestive enough to give an idea for other (more complex) structures. Curiously, despite the simplicity of the following statement, we could not trace it back in the literature.

\begin{num}\label{def-Lie} {\bf Deforming Lie algebras}:  If $(\mathfrak{g}, [\cdot, \cdot ])$ is a finite dimensional Lie algebra and
\[ H^{2}( \mathfrak{g}; \mathfrak{g})= 0 ,\]
where the cohomology groups are the Lie algebra cohomology with coefficients in the adjoint representation, then any smooth deformation
of $\mathfrak{g}$ is equivalent to the trivial one.
\end{num}

In other words, for any family $[\cdot, \cdot ]_{t}$ of Lie algebra structures on $\mathfrak{g}$, depending smoothly on the real parameter $t$ 
(in an open interval containing the origin), 
with $[\cdot, \cdot ]_{0}= [\cdot, \cdot ]$, there exists a smooth family of Lie algebra isomorphisms 
\[ h_{t}: \mathfrak{g}_{t}\rmap \mathfrak{g}, \ \ \text{where}\ \ \mathfrak{g}_{t}= (\mathfrak{g}, [\cdot, \cdot]_{t}),\]
defined for $t$ close enough to zero, and such that $h_0= Id$.

\begin{proof} We reduce the problem to the exactness
of certain cohomology \cite{NR} classes, which, in turn, will be implied by 
the perturbation lemma. Denote by 
$C^{*}= (C^{*}(\mathfrak{g}; \mathfrak{g}), b)$
the Chevalley-Eilenberg complex computing the Lie algebra cohomology of $\mathfrak{g}$ with coefficients in
$\mathfrak{g}$ (the adjoint representation). The derivatives of $[\cdot, \cdot]_{t}$ with respect to $t$ define closed cocycles
\[ c_{t}\in C^{2}(\mathfrak{g}_{t}; \mathfrak{g}_{t}) .\]
The $h_{t}$'s must satisfy $[h_t(v), h_t(w)]_0= h_t([v, w]_t)$. Applying $\frac{d}{dt}$, we must have:
\[ [a_t(v), w]_{t}+ [v, a_{t}(w)]_{t}= \frac{d}{dt} [v, w]_{t}+ a_{t}([v, w]_t) ,\]
where $a_{t}= h_{t}^{-1}\frac{d}{dt}h_t$. Conversely, starting with $a_t$ satisfying this 
equation, 
\[ h_{t}= exp(\int_{0}^{t} a_{s} ds) \in GL(\mathfrak{g}) \]
will produce the desired isomorphisms. Hence it suffices to show that 
the smooth family
of closed two-cocycles $-c_{t}$ can be transgressed to a smooth family of smooth one-cocycles $a_t$.
Let $\delta_{t}$ the Chevalley-Eilenberg boundary corresponding to $[\cdot, \cdot ]_{t}$. Then
\[ C^{*}(\mathfrak{g}_{t}; \mathfrak{g}_{t})= (C^{*}(\mathfrak{g}; \mathfrak{g}), d+ \delta_{t}) ,\]
where $\delta_t= \delta_t- \delta$. 
By hypothesis, we find $h, g$ defined on $C^{2}$ and $C^{3}$ 
contracting  $C^{*}$ in degree $2$. Since $\delta_0= 0$, $1- \delta_t h$ will be invertible
for $t$ small enough (use e.g. a norm on $\mathfrak{g}$), and we get contracting homotopies $H_t$ (in degree 2)
by Corollary \ref{main-cor}. Due to the explicit formula for $H_t$, and the fact that 
$(1- \delta_t h)^{-1}$ is smooth with respect to $t$, $a_t= H_{t}(c_t)$ will have the desired properties.
\end{proof}

\begin{num} {\bf Corollary}: Let $\mathfrak{g}$ be a vector bundle over a manifold $M$, with
Lie algebra structures on the fibers $\mathfrak{g}_{x}$, depending smoothly on $x\in M$.
If at $x_0\in M$ we have $H^{2}(\mathfrak{g}_{x_0}; \mathfrak{g}_{x_0})= 0$, then, around $x_0$, $\mathfrak{g}$
is isomorphic with the trivial bundle of Lie algebras with (Lie algebra) fiber $\mathfrak{g}_{x_0}$.
\end{num}

\begin{num} {\bf Corollary}: Any smooth deformation of a compact simply connected Lie group $G$ is equivalent to the trivial one.  
\end{num}

Here, a family $G_{t}$ of Lie groups with $G_0= G$ is smooth if they fit into a bundle of Lie groups
\cite{DoLa}, i.e. if $\mathcal{G}= \cup_t G_t$
can be given a manifold structure with the property that the projection into $t$ is a submersion, and the
group structure maps are smooth at the level of $\mathcal{G}$. 

\begin{proof} Since $G$ is compact and the foliation of $\mathcal{G}$ by the $G_t$'s has no holonomy,
Reeb stability allows us to assume that $G_t= G$, as manifolds, for all $t$ close enough to zero.
Passing to Lie algebras, we obtain a deformation of $\mathfrak{g}$. Since $\mathfrak{g}$ is semi-simple, 
$H^2(\mathfrak{g}; \mathfrak{g})$ must vanish, hence the Lie algebra deformation must be trivial. Since all $G_t$'s are simply connected,
we can integrate the resulting Lie algebra isomorphisms to a family $G_t\rmap G$ of Lie group isomorphisms.
\end{proof}

\begin{num} {\bf Remark: }\label{remark-Lie}\rm Regarding rigidity of Lie algebras, there is yet another way to go: there is a natural map $C$ (conjugation) which associates to $\phi\in GL(\mathfrak{g})$, the conjugation of $[\cdot, \cdot]$ with respect to $\phi$ viewed as an element in $C^{2}(\mathfrak{g}; \mathfrak{g})$, and then a map $J$ (Jacobi) which associates to $c\in C^{2}(\mathfrak{g}; \mathfrak{g})$, the expression that appears in the Jacobi identity for $c$, viewed as an element in $C^{3}(\mathfrak{g}; \mathfrak{g})$. One recovers $\delta: C^{1}(\mathfrak{g}; \mathfrak{g})\rmap C^{2}(\mathfrak{g}; \mathfrak{g})$ as the differential of $C$ at the identity matrix, and $\delta: C^{2}(\mathfrak{g}; \mathfrak{g})\rmap C^{3}(\mathfrak{g}; \mathfrak{g})$ as the differential of $J$ at the given bracket. A version of the inverse (or implicit) function theorem for finite dimensional manifolds (see e.g. Theorem on LG 3.22 in \cite{Se}) implies the following result which is closely related to Theorem \ref{def-Lie} (see also the next remark).

\begin{num} {\bf Corollary}\label{corollary-Lie}: If $H^{2}(\mathfrak{g}; \mathfrak{g})= 0$, then any Lie bracket on $\mathfrak{g}$ which is sufficiently closed to the original bracket $[\cdot, \cdot]$ of $\mathfrak{g}$, is conjugated to $[\cdot, \cdot]$.
\end{num}
\end{num}

\begin{num} \label{def-Lie-rk}{\bf Remark:}\rm \ We view \ref{def-Lie} (and its proof) as a ``Moser trick''-argument for Lie algebras, while \ref{remark-Lie} and
\ref{corollary-Lie} as an ``inverse-function theorem''-argument. Due to the simplicity of the structure we looked at here (Lie algebras), the proofs of \ref{def-Lie} and \ref{corollary-Lie} are closely related, and this is actually a point where the role of the perturbation lemma (in \ref{def-Lie}) is very close to that of Newton's method (in implicit function theorems). However, for more complex structures, the two type of arguments become quite different. On one hand, the ``inverse function theorem''-arguments become quite analytical, culminating with rigidity results proven by Nash-Moser techniques \cite{Ham} (e.g. used in \cite{Ham2} for foliations and in \cite{CrFe} for Poisson manifolds). On the other hand, the ``Moser trick''-arguments tend to give more geometric approaches. Actually, finding a more geometric proof of Conn's rigidity theorem for Poisson manifolds \cite{CrFe} was one of the reasons for looking at the relevance of the perturbation lemma to deformation problems. Although a geometric proof of Conn's result (based on Moser's trick and averaging) is possible \cite{CrFe2}, it seems to us that the most fruitful approach would be a more systematic study of graded Lie algebras and their Kuranishi spaces in the category of {\it tame Frechet spaces} (that are central to Nash-Moser techniques \cite{Ham}), similar to the known theory in the analytic category- see \cite{GM} and the references therein. (The need of tame Frechet spaces comes from the fact that the complexes controlling deformations are usually made of spaces of smooth sections of vector bundles, and the differential decreases the degree of differentiability).
Perturbation techniques would be encountered along the way (but, of course, the depth of such a theory would come from the use of smoothing operators).
\end{num}

\section{Deforming topological algebras}

In this section we point out the relevance of perturbation methods to deformations of topological algebras. 
We first recall the algebraic version of deformations \cite{Ger}, known as formal deformations. Given 
an algebra $(A, \cdot)$, a formal deformation of $A$ is a (associative) multiplication
\[ a\star b= ab + \sum_{k= 1}^{\infty} c_k(a, b) t^k \]
depending on the formal variable $t$. (For the notion of equivalence of such deformations, see \cite{Ger}). Such (equivalence classes of) formal deformations are ``controlled'' by $H^*(A; A)$, the Hocschild cohomology of
$A$ with coefficients in $A$. Recall that, for any $A$-bimodule $X$, $H^*(A; X)$ is computed by the complex
\[ C^k(A; X)= Hom(A^{\otimes k}, X),\]
with the boundary map $b: C^k(A; X)\rmap C^{k+1}(A; X)$,
\begin{eqnarray}
b(\alpha)(a_0, \ldots, a_k) & = & a_0 \alpha(a_1, \ldots, a_k)  \nonumber \\
                            & + & \sum_{i=1}^{k} (-1)^i \alpha(a_0, \ldots , a_{i-1}a_{i}, \ldots, a_k)+ \alpha(a_0, \ldots, a_{k-1})a_k.\nonumber
\end{eqnarray}
When $X= A$, the Hochschild complex $C^*(A; A)$ comes equipped with the Gerstenhaber bracket, 
\[ [\cdot, \cdot]: C^p(A; A)\times C^q(A; A)\rmap C^{p+q-1}(A; A), [\alpha, \beta]= \alpha\circ \beta- (-1)^{(p-1)(q-1)}\beta\circ \alpha ,\]
where 
\[ (\alpha\circ \beta)(a_1, \ldots , a_{p+q-1}) = \sum_{i= 0}^{p-1}(-1)^{(q-1)i} f(a_1, \ldots , g(a_{i+1}, \ldots , a_{i+q}),  \ldots a_{p+q-1}).\]
Moreover, $[\cdot, \cdot]$ passes to cohomology, and $(H^*(A; A), [\cdot, \cdot])$ is an algebraic analogue for the
space of multivector fields on a manifold $M$, endowed with the Nijenhuis-Schouten bracket (the two are the same if $A$ is the
algebra of smooth functions on $M$, and one imposes a certain continuity condition on the cocycles).

A similar discussion is possible for topological algebras. As in \cite{Tay} (to which we refer for more details also on the cohomology of topological algebras), we first fix a topological tensor product $\tilde{\otimes}$. The choice depends on the type of applications one has in mind, and it encodes the type of continuity one requires on the multiplication. For instance, joint continuity of the multiplication $m(a, b)= ab$ corresponds to the projective tensor product $\hat{\otimes}$ in the sense that $m$ is continuous if and only if it extends to a continuous map $m: A\hat{\otimes}A\rmap A$. Similarly, separate continuity corresponds to the inductive tensor product, etc. Having understood that we fix a topological tensor product $\tilde{\otimes}$, a topological algebra $A$ is an algebra endowed with a locally convex topology so that the multiplication $m$ of $A$ is continuous in the sense that it continuously extends to (or is defined as)  
\[ m: A\tilde{\otimes}A\rmap A .\]
The continuous cohomology, $H^{*}_{cont}(A; A)$, is defined as above, replacing $\otimes$ with $\tilde{\otimes}$. A formal deformation 
$\star$ of $A$ is called continuous if each of the coefficients $c_k$ are continuous (i.e. extend to $A\tilde{\otimes}A$). 
As in \cite{Ger}, if $H^{2}_{cont}(A; A)= 0$, then any continuous formal deformation of $A$ is equivalent to the trivial one.

In the topological world however, the deformations one would like to consider are families $\odot_{t}$ of (associative) multiplications so that
each $\odot_t$ is continuous, one recovers the original multiplication at $t= 0$, and $\odot_{t}\in Hom(A\tilde{\otimes}A, A)$ is of class
$C^k$ on $t$, where $0\leq k\leq \infty$ is fixed. We call such families {\it continuous deformations of $A$} which are of class $C^k$ on $t$.

\begin{num}{\bf Remark (on Poisson structures)}:\rm \ It is well-known that formal deformations of the algebra of smooth functions on a manifold $M$ are related to Poisson structures on $M$ (i.e. bivectors $\pi$ on $M$ satisfying $[\pi, \pi]= 0$, where $[\cdot, \cdot]$ is the Nijenhuis-Schouten bracket \cite{Vais}).  Motivated by the fact that $(H^*(A; A), [\cdot, \cdot])$ is the analogue of the space of multivector fields together with the Nijenhuis-Schouten bracket, a (non-commutative) Poisson structure on $A$ is defined \cite{Xu} by an element $\pi\in H^2(A; A)$ satisfying $[\pi, \pi]= 0$. Any formal deformation induces a Poisson structure, namely $\pi= [c_1]$. Indeed,
the associativity equation for $\star$ gives us certain equations corresponding to the powers of $t$. 
The coefficient of $t$ translates into $b(c_1)= 0$, while the one of $t^2$ translates into $[c_1, c_1]= b(c_2)$.

Similarly, continuous Poisson structures (i.e. with $\pi\in H^{2}_{cont}(A; A)$) are related to continuous deformation $\odot_{t}$.
More precisely, if $\odot_{t}$ is of class $C^1$ on $t$, then \[  c_1(a, b)= \frac{d}{dt}|_{t=0} a\odot_{t} b , \ \ (a, b\in A)\]
is a Hochschild cocycle, and $\pi= [c_1]\in H^2(A; A)$ is a candidate for a Poisson structure on $A$. If $\odot_{t}$ is of class $C^2$ on $t$, 
then $\pi$ is indeed a Poisson structure. To prove these assertions, one writes the associativity of the products $m_t(a, b)= a\odot_{t} b$ in terms of the Gerstenhaber bracket as
\[ [m_t, m_t]= 0,\]
and one uses that the Hochschild boundary is $b(\alpha)= [\alpha, m]$ ($m= m_0$ is the original multiplication). Taking derivatives with respect to
$t$ in the previous equation, one gets $[\frac{d}{dt}{m}_{t}, m_t]= 0$ hence $c_1= \frac{d}{dt}|_{t= 0}m_{t}$ is a cocycle. Taking derivatives at $t= 0$ in the last equation, we obtain $[c_1, c_1]= [c_2, m]$ hence exact, where
\[ c_2(a, b)= \frac{d^2}{dt^2}|_{t=0} a\odot_{t} b .\]
\end{num}

Exactly with the same proof as for Lie algebras (and using the fact that for Banach algebras small enough is small, cf. Remark \ref{in-rem} (iv)), we deduce:

\begin{num}  
\label{interesting}{\bf Proposition}: 
Let $A$ be a Banach algebra $A$ with the property that there are continuous maps $g$ and $h$ 
\begin{equation}  
\xymatrix{
C^{1}_{cont}(A; A) \ar@<-1ex>[r]_-{b} &   C^{2}_{cont}(A; A) \ar[l]_-{h}\ar@<-1ex>[r]_-{b} & C^{3}_{cont}(A; A)\ar[l]_-{g} ,\ bh+ gb= 1.}
\end{equation}
Then any continuous deformation of $A$ which is of class $C^1$ on $t$, is equivalent to the trivial one. 
\end{num}

A particular case is that of Banach algebras of (continuous) cohomological dimension at most one (see \cite{Tay}).

\begin{num} {\bf Corollary}: If the Banach algebra $A$ has cohomological dimension at most $1$, then any continuous deformation
of $A$ which is of class $C^1$ in $t$ is equivalent to the trivial one.
\end{num}

\begin{proof} As in the algebraic case, the condition on $A$ is equivalent (see \cite{Tay}) to the existence of a  
projective resolution $R$ of $A$ of type
\[ R: 0\rmap R_1\rmap R_0\rmap A .\] 
Then $Hom_{A-A}(R, A)$ and $C^{*}_{cont}(A; A)$ will be homotopic equivalent via continuous maps and homotopies, and this clearly produces a contracting homotopy of $C^{*}_{cont}(A; A)$ in degree $2$.
\end{proof}

\begin{num}{\bf Remark:}\rm\
We have seen that, in the case of Banach algebras, the hypothesis on the existence of a (linear continuous)
contraction $h$ can be relaxed (see Corollary \ref{cor-pert-Banach}), i.e. one has a 
weaker condition: $H^2(A; A)= )$ and $Im(b)\subset C^3(A; A)$ closed. This should be compared with the results
of \cite{Jo} and \cite{RT} which prove rigidity under perturbations under such conditions. 
The approach in \cite{RT} is an ``inverse-function theorem''-type approach (the paper
starts with such a theorem in this context). In contrast, 
the type of arguments in \cite{Jo} is a ``homological perturbation'' one. This is implicit in the proofs there,
and it is interesting to point out that a particular case of our Corollary \ref{cor-pert-Banach} already appears
(implicitly) in \cite{Jo}. Note also that such ``relaxations'' on the contracting homotopy $h$ are no longer possible for more
general topological vector spaces (not even for Frechet vector spaces), and, in the ``inverse function theorem''-approach to rigidity, 
this translates into that fact that the inverse function theorems one uses has similar hypothesis (e.g., for Nash-Moser techniques, see \cite{Ham2}, Theorem 3.1.1 in \cite{Ham2}, and Proposition 2.1 in \cite{CrFe}).
\end{num}

\section{Deforming algebras of continuous functions}

In this section we look at deformations $\odot_{t}$ of the product on the algebra of continuous function on a compact topological space $M$.
For the start, to avoid discussions on topological tensor products, we say that 
$\odot_{t}$ is continuous if it can be extended to a continuous map $C(M\times M)\rmap C(M)$. Smoothness with respect to $t$ 
will refer to smoothness as family of maps between the Banach spaces $C(M\times M)$ and $C(M)$.

\begin{num}\label{def-alg-cont} {\bf Deforming algebras of continuous functions}: Let $M$ be a compact metric space, and let $(C(M), \cdot)$ be the algebra of continuous functions on $M$.  Then any continuous deformation $(C(M), \odot_{t})$ of
the algebra $(C(M), \cdot)$ which is of class $C^1$ in $t$ is equivalent to the trivial one, i.e. there exists a family
\[ h_t: (C(M), \cdot)\rmap (C(M), \odot_{t})\]
of continuous algebra isomorphisms, defined for $t$ close enough to the origin, with $h_0= Id$, and which is of class $C^1$ in $t$. 
\end{num}


\begin{proof} The proof is similar to that for Lie algebras, using this time
the (Hochschild) complex computing the (continuous) cohomology of $A$ with coefficients in $A$,
$C^*= C^{*}_{cont}(A; A)$ (see the previous section).  
Continuity is the same type of continuity as for the products, 
i.e. one declares that the (topological) tensor product of $C(M)$ and $C(N)$ equals to $C(M\times N)$ for compact metric spaces $M, N$.
Hence
\[ C^{k}= Hom_{cont}(C(M^k), C(M)) .\]
And, proceeding as in the case of  Lie algebras, one only has to show that $C^*$ admits a continuous contraction in degree $2$.  
As in the algebraic case, the Hochschild complex comes from a projective resolution of $A$ by  $A$-bimodules, after applying the functor $Hom_{A-A}(-, A)$ of continuous $A$-bimodule maps. The resolution under discussion is
the $b'$-resolution:
\[\ldots \stackrel{b'}{\rmap}  C(M\times M\times M)\stackrel{b'}{\rmap} C(M\times M)\stackrel{m}{\rmap} C(M) \]
where $m$ is the multiplication (restriction to the diagonal), and
\[ b'(f_0, \ldots , f_{n+1})= \sum_{i=0}^{n} (-1)^i(f_0, \ldots, f_{i-1}, f_{i}f_{i+1}, f_{i+2}, \ldots , f_{n+1}) .\]
To prove that $C^*$ admits a continuous contraction in degree $2$, it suffices to show that the $b'$-resolution admits a similar contraction, which is a map of $A$-bimodules. On the other hand, one can replace the $b'$ resolution by its normalization $N^*$, i.e. consisting on those continuous functions satisfying
\[ f(x_0, \ldots, x_n)= 0, \ \text{if}\ x_i= x_{i+1}\ \text{for\ some}\ 0\leq i\leq n-2 .\]
Indeed, the fact that a complex (associated to a simplicial module) and its normalization have the same cohomology can be proven by giving the explicit formulas for the maps and the homotopies (see the proof of Theorem 6.1. in \cite{ML}), and these maps are clearly continuous. Finally, using a metric $\rho$, we can write down an explicit formula for a contracting homotopy $h$ for
\[ \ldots  \stackrel{b'}{\rmap} N^2 \stackrel{b'}{\rmap} N^1 \stackrel{b'}{\rmap} ker(m) .\]
We give the first three components of $h$ (which is a bit more that what we need), and a general formula can be easily guessed. For $f\in Ker(m)\subset C(M^2)$, $g\in C(M^3)$, $k\in C(M^4)$:
\begin{eqnarray}
h(f)(x_0, x_1, x_2) & = & \frac{\rho(x_0, x_1)}{\rho(x_0, x_1)+ \rho(x_1, x_2)} f(x_0, x_2),\nonumber \\
h(g)(x_0, x_1, x_2, x_3) & = & \frac{\rho(x_0, x_1)\rho(x_1, x_2)}{\rho(x_0, x_1)\rho(x_1, x_2)+ \rho(x_2, x_3)}\big{(} g(x_0, x_1, x_3)\nonumber\\
   & - & \frac{\rho(x_0, x_1)}{\rho(x_0, x_1)+ \rho(x_1, x_2)}g(x_0, x_2, x_3)\big{)}\nonumber \\
h(k)(x_0, x_1, x_2, x_3, x_4) & = & \frac{\rho(x_0, x_1)\rho(x_1, x_2)\rho(x_2, x_3)}{\rho(x_0, x_1)\rho(x_1, x_2)\rho(x_2, x_3)+ \rho(x_3, x_4)} \big{[} k(x_0, x_1, x_2, x_4) \nonumber \\
        & - & \frac{\rho(x_0, x_1)\rho(x_1, x_2)}{\rho(x_0, x_1)\rho(x_1, x_2)+ \rho(x_2, x_3)}\big{(}k(x_0, x_1, x_4, x_4)\nonumber\\
   & - & \frac{\rho(x_0, x_1)}{\rho(x_0, x_1)+ \rho(x_1, x_2)} k(x_0, x_2, x_4, x_4)\big{)}\big{]} \nonumber
\end{eqnarray}
\end{proof}

\begin{num}\label{def-Ban-rk}{\bf Remark:}\rm \ The previous result is similar to the fact that the cyclic homology of $C(M)$ is un-interesting (in contrast with the algebra of smooth functions on a manifold $M$,  when one obtains the usual DeRham cohomology \cite{Connes}). Actually, the previous proof 
shows that $C(M)$ is ``cohomologically un-interesting'' and it does imply the triviality of the cyclic homology. 

Incidentally, let us also point out that the vanishing of the Hochschild cohomology 
and the rigidity theorem of Section 3 of \cite{Ger} adapted to the topological context (see also the next section) also implies the following:
\end{num}

\begin{num} {\bf Corollary}: Given a compact metric space $M$, any formal deformation of $C(M)$:
\[ f\star_{t}g= fg+ c_1(f, g)t+ c_{2}(f, g)t^2+\ldots \ ,\]
with the property that each of the coefficients $c_k$ is continuous, is equivalent to the trivial one.
\end{num}

\end{document}